\documentclass[10pt]{article}
\font\smallit=cmti10
\font\smalltt=cmtt10

\usepackage{amssymb,latexsym,amsmath,epsfig,amsthm, hyperref, cleveref} 

\makeatletter

\renewcommand\section{\@startsection {section}{1}{\z@}
{-30pt \@plus -1ex \@minus -.2ex}
{2.3ex \@plus.2ex}
{\normalfont\normalsize\bfseries\boldmath}}

\renewcommand\subsection{\@startsection{subsection}{2}{\z@}
{-3.25ex\@plus -1ex \@minus -.2ex}
{1.5ex \@plus .2ex}
{\normalfont\normalsize\bfseries\boldmath}}

\renewcommand{\@seccntformat}[1]{\csname the#1\endcsname. }

\DeclareMathOperator{\ord}{ord}

\newcommand{\Z}{\mathbb{Z}}
\newcommand{\Q}{\mathbb{Q}}

\newcommand{\C}{\mathbb{C}}

\renewcommand{\aa}{\alpha}

\newcommand{\ps}[1]{\left(#1\right)}
\newcommand{\floor}[1]{\left\lfloor #1 \right\rfloor}

\newcommand{\poc}[1]{(#1; #1)_{\infty}}

\makeatother

\newtheorem{theorem}{Theorem}
\newtheorem{lemma}{Lemma}

\theoremstyle{definition}

\newtheorem*{remark*}{Remark}
\newtheorem*{example*}{Example}

\begin{document}

\begin{center}
\uppercase{\bf Congruences For Fractional Partition Functions}
\vskip 20pt
{\bf Yunseo Choi}\\
{\smallit Phillips Exeter Academy, Exeter, NH, USA}\\
{\tt ychoi@exeter.edu}\\
\end{center}
\vskip 20pt

\centerline{\smallit Received: 8/11/19, Revised: 10/17/20, Accepted: 1/9/21, Published: 2/1/21}
\vskip 30pt 

\centerline{\bf Abstract}

\noindent
The coefficients of the generating function $(q;q)^\alpha_\infty$ produce $p_\alpha(n)$ for $\alpha \in \mathbb{Q}$. In particular, when $\alpha = -1$, the partition function is obtained. Recently, Chan and Wang studied congruences for $p_\alpha(n)$ and gave several infinite families of congruences of the form $p_\alpha(\ell n + c) \equiv 0 \pmod{\ell}$ for primes $\ell$ and integers $c$. Expanding upon their work, given adequate $\alpha$, we use the lacunarity of the powers of the Dedekind-eta function to raise the modulus of Chan and Wang's congruences to higher powers of $\ell$. In addition, we generate new infinite classes of congruences through the multiplicative properties of the coefficients of Hecke eigenforms. This allows us to prove new families of congruences such as: $p_{-\frac{1}{8}}(7^2n+5)\equiv 0 \pmod{7^2}$. 

\pagestyle{myheadings}
\markright{\smalltt INTEGERS: 20 (2020)\hfill}
\thispagestyle{empty}
\baselineskip=12.875pt
\vskip 30pt 

\section{Introduction}

A $\textit{partition}$ of a non-negative integer $n$ is a non-increasing sequence of positive integers that sum to $n$. Per usual, let $p(n)$ denote the number of distinct ways to partition $n$. Euler discovered the generating function of the partition function to be: 
\begin{equation*}
    P(q):= \sum_{n=0}^{\infty} p(n)q^n = \frac{1}{\poc{q}},
\end{equation*}
where $\poc{q} := \prod_{n=1}^{\infty} (1-q^{n})$ is the $q$-Pochhammer symbol, defined for $\lvert q \rvert < 1$.

Ramanujan observed and proved congruences in $p(n)$ for $n$ in special arithmetic progressions.
\begin{equation*}
    p( 5n + 4 ) \equiv 0 \pmod{5},
\end{equation*}
\begin{equation*}
    p( 7n + 5 ) \equiv 0 \pmod{7},
\end{equation*}
\begin{equation*}
    p( 11n + 6 ) \equiv 0 \pmod{11}. 
\end{equation*}

In addition, Ramanujan conjectured that for all powers of $\ell \in \{5, 7, 11\}$, there exists a class of congruences in which the common difference of the arithmetic progression and the modulus share the same power of $\ell$. His conjecture was proven to be false when Chowla and Gupta \cite{web} discovered $7^3$ to be a counterexample. Nonetheless, a slight modification of the conjecture was proven to hold true by Atkin \cite{atkin} and Watson \cite{watson}: for any $k \in \Z^+$ and a prime $\ell \in \{5, 7, 11\}$, when $r_{\ell,k} \equiv 1/24 \pmod{\ell^k},$ we have for all $n$ that
\begin{align*}
    p(5^kn + r_{5,k}) &\equiv 0 \pmod{5^k},\\
    p(7^kn + r_{7,k}) &\equiv 0 \pmod{7^{\floor{k/2}+1}},\\
    p(11^kn + r_{11,k}) &\equiv 0 \pmod{11^k}.
\end{align*}

When the condition that the common difference of the arithmetic sequence and the modulus have to be the powers of the same prime is relaxed, many more congruences are present. In fact, Ono and Ahlgren \cite{ahlono} proved that for all integers $L$ co-prime to 6, there exist $A, B \in \Z$ such that for all $n$, $p(An+B) \equiv 0 \pmod{L}$. 

The continued search for congruence relations in the partition function led to the search of congruence relations in fractional partition functions. The fractional partition function is the generating function of the usual partition function raised to the power of $-\alpha \in \Q$. Throughout this paper, we let $\alpha=\frac{a}{b}$ where $\alpha$ is a fraction written in lowest terms with a positive denominator. Let
\begin{equation*}\label{fracp}
   P_{\alpha}(q):= {\poc{q}^{\alpha}}
   := \sum_{n=0}^{\infty} p_{\alpha}(n)q^n.
\end{equation*}

We set $p_\aa(n) := 0$ for $n < 0.$ Unlike $p(n)$ that are integral, $p_{\alpha}(n)$ is a non-integral rational number for most choices of $n$ and $\alpha$. Chan and Wang \cite{chan} addressed this issue in the context of congruences (Theorem 1.1 of \cite{chan}) by showing that that $p_\aa(n)$ are $\ell$-integral for any prime $\ell \nmid b$.

In addition, Chan and Wang (Theorem 1.2 of \cite{chan}) displayed infinite families of congruences for fractional partition functions, making use of the previously-known, explicit expressions of the coefficients of $\poc{q}^{d}$ for $d \in \{1, 3, 4, 6, 8, 10, 14, 26\}$.
\begin{theorem}
    \label{chan_wang_congruences}
    (Cf. [4, Theorem 1.2]) Supposed that $\alpha \in \mathbb{Q}$ and $d, r \in \mathbb{Z}$ are given. Let $\ell$ denote a prime such that $\ell \mid a-db$. If $d, r,$ and $\ell$ satisfy one of the following conditions: 
\begin{enumerate}
    \item $d=1$ and $(\frac{24r+1}{\ell})= -1$;
    \item $d=3$ and $(\frac{8r+1}{\ell}) \ne 1$;
    \item $d \in \{4, 8, 14\}$, $\ell \equiv 5 \pmod{6}$ and $\ell \mid 24r+d$;
    \item $d \in \{6, 10\}$, $\ell \geq 7$, $\ell \equiv 3 \pmod{4}$ and $\ell \mid 24r + d$;
    \item $d = 26$, $\ell \equiv 11 \pmod{12}$ and $\ell \mid 24r+d$, 
\end{enumerate}
then, for all $n$, we have that $p_{\alpha}(\ell n + r) \equiv 0 \pmod{\ell}$. 
\end{theorem}

Notice that for each $d \in \{4, 6, 8, 10, 14, 26\}$, there are conditions imposed on $\ell$ independent of the choices of $\alpha$ and $r$. For example, when $d \in \{6, 10\}$, it is required that $\ell \geq 7$ and that $\ell \equiv 3 \pmod{4}$. For each $d$, we define a prime $\ell$ to be \textit{$d$-satisfactory} if $\ell$ satisfies such exact conditions, except that we additionally exclude 5 from the list of \textit{$14-$satisfactory} primes and 11 from the list of \textit{$26-$satisfactory} primes.   

It is natural to ask about the significance of the list of $d$ in Chan and Wang's theorem. This brings us to a result by Serre \cite{lacunary} on Dedekind eta-functions, defined as ${\eta(\tau) := q^{1/24} \poc{q}}$ for $q := e^{2\pi i\tau}.$ Recall that a Fourier expansion $\sum_{n = 0}^{\infty} a(n)q^n$ is \textit{lacunary} if 
\[ \lim_{N \to \infty} \frac{\#\{n \le N : a(n) = 0\}}{N} = 1. \] In 1985, Serre \cite{lacunary} proved that $\eta(\tau)^{d}$ is lacunary for $d \in 2\mathbb{Z}$ if and only if $d \in \{2, 4, 6, 8, 10, 14, 26\}$. In addition, Serre provided explicit ways of writing such lacunary $\eta$ powers as linear combinations of Hecke eigenforms.  

In \Cref{first_theorem}, we make use of Serre's results on the lacunarity of $\eta$ powers and raise the power of $\ell$ in the modulus of Chan and Wang's congruences to $ord_{\ell}(\alpha-d)$. In other words, given our choice of $\alpha$, the power of $\ell$ in the modulus can be arbitrarily high. 

\begin{theorem}
\label{first_theorem}
    For $d \in \{4, 6, 8, 10, 14, 26\}$, let $\ell$ be a \textit{$d$-satisfactory} prime. If $r$ satisfies $ord_{\ell}(\frac{24}{\gcd(d,24)} r+ \frac{d}{\gcd(d,24)}) =  1$, then, we have for all $n$ that
    \begin{equation*}
        p_{\alpha}(\ell^2 n + r) \equiv 0 \pmod{\ell^{ord_{\ell}(\alpha-d)}}. 
    \end{equation*}
\end{theorem}

\begin{remark*}
    Although $5$ and $11$ were removed from the list of $14$ and $26-$\textit{satisfactory} primes, a modified statement of \Cref{first_theorem}--that is, the power of $\ell$ in the modulus is not $ord_\ell(\alpha-d)$, but instead is $\ord_{\ell} (\alpha-d)-1$ and $\ord_{\ell} (\alpha-d)-2$, respectively--holds true for such choices of $d$ and $\ell$ (See \Cref{eta_function}).
\end{remark*}

\begin{example*}
We demonstrate an example and show that for certain choices of $\alpha, r, d,$ and $\ell$, the power of $\ell$ in the modulus given by \Cref{first_theorem} is sharp. Let $\ell=7$ and $d=6$. $\ell$ is \textit{$6-$satisfactory} because $\ell \geq 7$ and $\ell \equiv 3 \pmod{4}$. In addition, we let $r=5$ as $ord_{7}(4 \cdot 5 + 1) = 1$. Now, let $\alpha = -\frac{1}{8}$. Since $\ord_{7}({-\frac{1}{8}}-6) = 2$, we conclude from \Cref{first_theorem} that 
\begin{align*}
    p_{-\frac{1}{8}}(7^2 n + 5) &\equiv 0 \pmod{7^2}.
\end{align*}
The power of $7$ in the modulus given by \Cref{first_theorem} is sharp in this case because
\[ p_{-\frac{1}{8}}(7^2 \cdot 0 + 5) = p_{-\frac{1}{8}}(5) \equiv \frac{55615}{262144} \not\equiv 0 \pmod{7^3}. \]
\end{example*}

It is also conspicuous that while many integers in Chan and Wang's list and Serre's list coincide, $d = 2$ is missing from Chan and Wang's list. We cover this case in \Cref{second_theorem} by showing that a slightly weaker statement of \Cref{first_theorem} holds true for $d=2$. We define a prime $\ell$ to be \textit{$2$-satisfactory} if $\ell \not\equiv 1 \pmod{12}$.

\begin{theorem}
    \label{second_theorem}
    For $d = 2$, let $\ell$ be a $2$-\textit{satisfactory} prime. If $r$ satisfies $ord_{\ell}(12r + 1) =  1$, then, we have for all $n$ that
    \begin{equation*}
     p_{\alpha}(\ell^2 n + r) \equiv 0 \pmod{\ell^{ord_{\ell}(\alpha-2)-1}}.  
    \end{equation*}
\end{theorem}

\begin{example*}
We once again give an example and show that for certain choices of $\alpha,$ $r,$ and $\ell,$ the power of $\ell$ in the modulus given by \Cref{second_theorem} is sharp. Let $\ell = 5$, a \textit{$2$-satisfactory} prime as $\ell \not\equiv 1 \pmod{12}$. Let $\alpha=\frac{1}{13}$. Since $\ord_{5}({\frac{1}{13}}-2) = 2$, it follows from \Cref{second_theorem} that 
\begin{align*}
    p_{\frac{1}{13}}(5^2 n + 7) &\equiv 0 \pmod{5^1}.
\end{align*}
The power of $5$ in the modulus given by \Cref{second_theorem} is sharp in this case because
\[ p_{\frac{1}{13}}(5^2 \cdot 0 + 7) = p_{\frac{1}{13}}(7) \equiv -\frac{ 3395395}{62748517} \not\equiv 0 \pmod{5^2}. \]

\end{example*}

\Cref{first_theorem} and \ref{second_theorem} rely heavily on the lacunarity of the corresponding $\eta$ powers (See \Cref{section:proofs}). For $d=2$, however, adequate choices of arithmetic progressions along the coefficients of $\eta(12\tau)^2$ produce sequences with elements that are not uniformly 0, but are nonetheless the multiples of the same prime power. This leads us to our final theorem. 

\begin{theorem}
\label{third_theorem}
    For $d=2$, fix a prime $\ell$ and $v \in \Z^{+}$. Then, there exists a finite $w \in \Z^{+}$ such that when $ord_{\ell} (\alpha-2) = v+w$ and  $ord_{\ell}(12r + 1)  = w$, we have for all $n$ that
    \begin{equation*}
        p_{\alpha} (\ell^{w+1}n+ r) \equiv 0 \pmod{\ell^v}. 
    \end{equation*}
\end{theorem}

\begin{remark*}
The significance of \Cref{third_theorem} is that we may drop the \textit{$2$-satisfactory} condition. If $\ell$ is \textit{$2$-satisfactory}, \Cref{second_theorem} and \ref{third_theorem} give the same congruences. 
\end{remark*}

\begin{example*}
We provide an example that is not covered by \Cref{second_theorem} by choosing an $\ell$ that is not \textit{$2$-satisfactory}. $\ell=13$ is one such prime, and we let $v=1$. Then, we show that $w=12$ is a valid choice of $w$ (See \Cref{third_proof_lemma}). Computation on \verb Mathematica  shows that $a_{2}(1) = 1$ and $a_{2}(13)=-2$. Now, setting $\ell=13$ in Equation \eqref{eq:3.11} gives that $a_{2}(13^{k})=(-1)^{k+1}(k+1)$ for $k \in \mathbb{Z}^{+}$. In particular, we have that $a_{2}(13^{12}) \equiv 0 \pmod{13}$. We let $r= \frac{11\cdot13^{12}-1}{12}$, since $ord_{13}(12 \cdot \frac{11\cdot13^{12}-1}{12}+1) = 12$. In addition, note that for $\alpha=\frac{a}{b}$ such that $a=1$ and $b=\frac{13^{13}+1}{2}$, $ord_{\ell}(\alpha -2) = ord_{\ell}(\frac{2}{13^{13}+1} -2) = ord_{\ell}(\frac{-2 \cdot 13^{13}}{13^{13}+1}) = 13$. Thus, for such $\alpha$, \Cref{third_theorem} gives for all $n$ that
\[ p_{\alpha}(13^{13} \cdot n + \frac{11\cdot13^{12}-1}{12}) \equiv 0 \pmod{13^{1}}. \]
\end{example*}

\section{Preliminaries}
\label{preliminaries}
\subsection{Modular Forms} \label{modular_forms} These facts are well-known and can be found in any standard text, such as \cite{web}. First, we define the Eisenstein series that describes modular forms. To do so, we define the divisor function $\sigma_{k-1}(n)$ for positive integers $k$:
\begin{equation*}
    \sigma_{k-1}(n) := \sum_{1 \leq d \mid n} d^{k-1}. 
\end{equation*}

Now, recall that all modular forms of $SL_2(\Z)$ are generated by $E_4(\tau)$ and $E_6(\tau)$ where: 
\begin{equation*}
    E_4(\tau) = 1 + 240 \sum\limits_{n=1}^{\infty} \sigma_3(n)q^n     \text{  and}
\end{equation*}
\begin{equation*}
    E_6(\tau) = 1 - 504 \sum\limits_{n=1}^{\infty} \sigma_5(n)q^n.
\end{equation*}

Next, we define the congruence subgroup of $SL_2(\Z)$  of level $N$, denoted by $\Gamma_0(N)$. 
\begin{equation*}
    \Gamma_0(N) = \Big\{\left(\begin{array}{cc} a & b\\ c & d \end{array}\right) \in SL_2(\Z) : c \equiv 0 \pmod{N}\Big\}.
\end{equation*}

In addition, we let $M_k(\Gamma_0(N))$ refer to the complex vector space of modular forms of weight $k$ with respect to $\Gamma_0(N)$. If $\chi$ is a Dirichlet character modulo $N$, we say that a modular function $f(\tau) \in M_k(\Gamma_0(N))$ has a \textit{Nebentypus character} $\chi$ if for all $\tau \in \mathbb H$ and for all $\left(\begin{array}{cc} a & b\\ c & d \end{array}\right) \in \Gamma_0(N)$,
\begin{equation*}
    f(\frac{a\tau+b}{c\tau+d}) = \chi(d) (c\tau+d)^k f(\tau).
\end{equation*}
The space formed by such modular forms is referred to as $M_k(\Gamma_0(N), \chi)$. Additionally, we note that the $m$th Hecke operator for $m \in \Z^+$, $T_{m, k, \chi}$, is an endomorphism on $M_k.$ Its action on a Fourier expansion $f(\tau) = \sum_{n = 0}^{\infty} a(n)q^n$ is illustrated by the formula:
\begin{equation*}
    f(\tau) \mid T_{m,k, \chi} = \sum_{n=0}^{\infty} \ps{\sum_{\delta\mid (m,n)}
    \chi(\delta)\delta^{k-1}a(mn/\delta^2)}q^n.
\end{equation*}
When $m = \ell$ is a prime, the expression reduces to \[ f(\tau)\mid T_{\ell, k, \chi} = \sum_{n=0}^{\infty} \ps{a(\ell n) + \chi(\ell)\ell^{k-1}a(n/\ell)}q^n, \]
where $a(\frac{n}{\ell})=0$ for $\ell \nmid n$.
Recall that a modular form $f(\tau) \in M_k(\Gamma_0(N), \chi)$ is a \textit{Hecke eigenform} if it is an eigenvector of $T_{m, k, \chi}$ for all $m \ge 1$, i.e. if there exist a $\lambda(m) \in \C$ such that \[ f(\tau)\mid T_{m,k} = \lambda(m)f(\tau). \] 

 In particular, if $a(1)=1$, then we consider $f(\tau)$ to be \textit{normalized}. This definition naturally leads us to the following lemma. The proof of this lemma follows immediately from the definitions. 

\begin{lemma}
\label{first_lemma}
Suppose that $f (\tau) = \sum_{n = 0}^{\infty} a(n)q^n \in M_k(\Gamma_0(N), \chi)$ is a normalized cuspidal Hecke eigenform. Then, it follows that
\begin{equation*}
    a(n)a(\ell) = a(n\ell) + \chi(\ell)\ell^{k-1}a(\frac{n}{\ell}).
\end{equation*}
\end{lemma}

\subsection{On the Powers of the Dedekind Eta Function} \label{eta_function}

The Dedekind eta function is defined as ${\eta(\tau) := q^{1/24} \poc{q}}$ for $q := e^{2\pi i\tau}.$
It is known by Martin \cite{martin} that $\eta(\tau)^d$ for $d \in \{1, 2, 3, 4, 6, 8, 12, 24\}$ are Hecke eigenforms. In addition, Carney, Etropolski, and Pitman (Lemma 2.2 of \cite{eigenform}) characterized $\chi(d)$ for each $\eta(\tau)^d$. 

\begin{lemma}
\label{second_lemma}
    $\chi(d)$ for $\eta(\tau)^d$ for $d \in \mathbb{Z}$ is 
    \begin{equation*}
    \chi(d) := \begin{cases} 
      (\frac{(-1)^{\frac{d}{2}}}{\cdot}) & \text{if } d \in 2\Z \\
      (\frac{12}{\cdot}) & \text{if } d \not\in 2\Z \cup 3\Z \\
      (\frac{-4}{\cdot}) & d \in 3\Z \setminus 2\Z.
   \end{cases}
\end{equation*}
\end{lemma}

In 1985, Serre \cite{lacunary} proved that $\eta(\tau)^d$ for $d \in 2\mathbb{Z}$ is lacunary if and only if $d \in \{2, 4, 6, 8, 10, 14, 26\}$. Additionally, for each of such $d$, he presented explicit ways to write $\eta(\frac{24}{gcd(d, 24)}\tau)^d$ in linear combinations of Hecke eigenforms. The expression $\frac{24}{gcd(d,24)}$, multiplied to $\tau$, ensures that $\eta(\frac{24}{gcd(d, 24)}\tau)^d$ is an expression of integral powers of $q$. As the specifics of these formulae play an integral role in proving our results, we list the formulae. In addition, we note that throughout the paper, we denote $\eta(\frac{24}{gcd(d,24)}\tau)^d = \sum_{n = 0}^\infty a_{d}(n)q^{n}$. 

If $d \in \{2, 4, 6, 8, 12\}$, $\eta(\tau)^d$ are Hecke eigenforms themselves. For $d= 10$, $\eta(12\tau)^{10}$ can be written as a linear combination of two Hecke eigenforms, $E_4(12\tau)\eta(12\tau)^2 \pm 48 \eta(12\tau)^{10}$. We have 
\begin{equation} \label{eq:2.1}
    \begin{aligned}
        \eta^{10}(12\tau)  = \frac{1}{96} ((E_4(12\tau)\eta(12\tau)^2 + 48 \eta(12\tau)^{10}) \\ - (E_4(12\tau)\eta(12\tau)^2 - 48 \eta(12\tau)^{10}) ).
    \end{aligned}
\end{equation}

Note that because $10-$\textit{satisfactory} primes $\ell$ are co-prime with $96$, the factor of $\frac{1}{96}$ does not interfere with divisibility modulo $\ell$. 

Similarly, $\eta(12\tau)^{14}$ is a linear combination of two Hecke eigenforms, namely, $E_6(12\tau) \eta(12\tau)^2 \pm 360\sqrt{-3}\eta(12\tau)^{14}$. We have 
\begin{equation}
\begin{aligned} \label{eq:2.2}
\eta(12\tau)^{14} = & \frac{1}{720\sqrt{-3}} ((E_6(12\tau) \eta(12\tau)^{2} + 360\sqrt{-3}\eta(12\tau)^{14})\\
      & -(E_6(12\tau) \eta(12\tau)^2 - 360\sqrt{-3}\eta(12\tau)^{14})).\\
\end{aligned}
\end{equation}
 We remove $5$ from the list of $14-$\textit{satisfactory} primes, because the constant factor of $\frac{1}{720}$ divides out a factor of $5$ from the numerator. 
 
For $d= 26$, $\eta(12\tau)^{26}$ can be written as a linear sum of four Hecke eigenforms, specifically, $E_6^2(12\tau)\eta(12\tau)^2 + 9398592 \eta(12\tau)^{26} \pm 102960 \sqrt{-3}E_6(12\tau)\eta(12\tau)^{14}$ and $E_6^2(12\tau)\eta(12\tau)^2 - 6910272 \eta(12\tau)^{26} \pm 20592E_8(12\tau)\eta(12\tau)^{10}$. We have
\begin{equation}
\begin{aligned} \label{eq:2.3}
\eta(12\tau)^{26} = & \frac{1}{32617728} ((E_6^2(12\tau)\eta(12\tau)^2 + 9398592 \eta(12\tau)^{26} \\ & + 102960 \sqrt{-3}E_6(12\tau)\eta(12\tau)^{14}) 
       +(E_6^2(12\tau)\eta(12\tau)^2 \\ & + 9398592 \eta(12\tau)^{26} - 102960 \sqrt{-3}E_6(12\tau)\eta(12\tau)^{14}) \\
      & -(E_6^2(12\tau)\eta(12\tau)^2 - 6910272 \eta(12\tau)^{26} + 20592E_8(12\tau)\eta(12\tau)^{10}) \\
      & -(E_6^2(12\tau)\eta(12\tau)^2 - 6910272 \eta(12\tau)^{26} - 20592E_8(12\tau)\eta(12\tau)^{10})).
\end{aligned}
\end{equation}
For the same reason that we removed $5$ from the list of $14-$\textit{satisfactory} primes, we remove $11$ from the list of $26-$\textit{satisfactory} primes. 

\subsection{Preliminary Results}

We state two key results by Chan and Wang \cite{chan}. The first result (Theorem 1.1 of \cite{chan}) identifies the congruences that are meaningful to study. 

\begin{theorem}
When written in lowest terms, we have that
\begin{equation*}
    \text{denom}(p_{\alpha}(n))=b^n\prod_{p\mid b}p^{\ord_p(n!)}.
\end{equation*}
\end{theorem} 

In other words, $\text{denom}(p_{\alpha}(n))$ is $\ell$-integral for any prime $\ell \nmid b$. We thus conclude that for a given rational number $\alpha$, whenever $\gcd(\ell, b) = 1,$ congruences modulo $\ell$ and its powers are well-defined. 

The second result is a technical lemma (Lemma 2.1 of \cite{chan}) resulting from Frobenius endomorphism. This lemma allows us to move exponents through $q$-Pochhammer symbols, a crucial step in the proofs of our main results.
\begin{lemma}
\label{third_lemma}
Let $\ell$ be a prime such that $\ell \nmid b$ as usual. Then, for any $r \ge 1$, we have that 
\begin{align*}
    \poc{q}^{\ell^r\aa} \equiv \poc{q^\ell}^{\ell^{r - 1}\aa} \pmod{\ell^r}.
\end{align*}
\end{lemma}

\section{Proofs of the main results} \label{section:proofs}

\begin{proof}[Proof of \Cref{first_theorem}]

We work out the case of $d=4$. Similar conclusions can be made about $d= 6$ and $8$ by following the same steps. For simplicity, we write $v := \ord_\ell(\aa-4)$ such that $\aa - 4 = \ell^{v} u$ for some $u \in \Z_{(\ell)}.$ First, we relate $p_{\alpha}(n)$ to $\eta(6\tau)^4$ using the $q$-Pochhammer symbol. We have that
\begin{equation}
    \begin{aligned}
    \sum_{n = 0}^\infty p_{\aa}(n)q^{6n + 1} &= q\poc{q^{6}}^{\alpha} =  q\poc{q^{6}}^{\ell^{v} u + 4} \\
    &= q\poc{q^{6}}^{4}\poc{q^{6}}^{\ell^{v} u} = \eta(6\tau)^{4}\poc{q^{6}}^{\ell^{v} u}.
    \end{aligned}
\end{equation}

Now, applying \Cref{third_lemma}, we have that   
\begin{equation} \label{eq:3.2}
    \sum_{n = 0}^\infty p_{\aa}(n)q^{6n + 1} = \eta(6\tau)^{4}\poc{q^{6}}^{\ell^{v} u} \equiv \eta(6\tau)^{4}\poc{q^{6\ell}}^{\ell^{v-1} u}
   \pmod{\ell^v}. 
\end{equation}

Recall that $\eta(6\tau)^4 = \sum_{n = 0}^\infty a_{4}(n)q^{n}$, and let $r_0$ denote the smallest positive integer such that $6r_0+1 \equiv 0 \pmod{\ell}$. Extracting the terms of the form $q^{\ell n}$ from both sides of Equation \eqref{eq:3.2} and replacing $q^{\ell}$ with $q$, we arrive at 
\begin{equation} \label{eq:3.3}
    \sum_{n = 0}^\infty p_{\aa}(\ell n+ r_0)q^{6n + \frac{6r_0 +1}{\ell}} \equiv \sum_{n = 0}^\infty a_4(\ell n)q^{n}\cdot \poc{q^{6}}^{\ell^{v-1} u}
   \pmod{\ell^v}. 
\end{equation}

Since $\ell$ is $4$-satisfactory and because $6r_0+1 \equiv 0 \pmod{\ell}$, it follows from \Cref{chan_wang_congruences} that $p_{\alpha}(\ell n+r_0) \equiv 0 \pmod{\ell}$. This allows us to divide each side of Equation \eqref{eq:3.3} by $\ell$. We now have that
\begin{equation} \label{eq:3.4}
    \frac{1}{\ell} \cdot \sum_{n = 0}^\infty p_{\aa}(\ell n+ r_0)q^{6n + \frac{6r_0 +1}{\ell}} \equiv \frac{1}{\ell} \cdot \sum_{n = 0}^\infty a_4(\ell n)q^{n}\cdot \poc{q^{6}}^{\ell^{v-1} u}
   \pmod{\ell^{v-1}}. 
\end{equation}

We apply \Cref{third_lemma} again and deduce that
\begin{equation} \label{eq:3.5}
     \frac{1}{\ell} \cdot \sum_{n = 0}^\infty p_{\aa}(\ell n+ r_0)q^{6n + \frac{6r_0 +1}{\ell}} \equiv \frac{1}{\ell} \cdot \sum_{n = 0}^\infty a_4(\ell n)q^{n}\cdot \poc{q^{6\ell}}^{\ell^{v-2} u}
   \pmod{\ell^{v-1}}. 
\end{equation}

Multiply $\ell$ back on both sides of Equation \eqref{eq:3.5} to arrive at 
\begin{equation} \label{eq:3.6}
     \sum_{n = 0}^\infty p_{\aa}(\ell n+ r_{0})q^{6n + \frac{6r_{0} +1}{\ell}} \equiv \sum_{n = 0}^\infty a_4(\ell n)q^{n}\cdot \poc{q^{6\ell}}^{\ell^{v-2} u}
   \pmod{\ell^{v}}. 
\end{equation}

Recall that $\frac{\eta(6\tau)^4}{q}$ is expression of $q^{6}$. As a result, $a_{4}(\ell) = 0$ for $\ell \equiv 5 \pmod{6}$. In addition, because $\eta(6\tau)^4$ is a normalized Hecke eigenform, it follows from \Cref{first_lemma} that it has multiplicative coefficients for co-prime indices, i.e., for any $k \in \mathbb{Z}_{(\ell)}$, we have that $a_{4}(\ell k) = 0$.

Finally, we extract the terms of the form $q^{\ell n + \frac{6r+1}{\ell}}$ from each side of Equation \eqref{eq:3.6}. Because $ord_{\ell} (6r+1) = 1$, $\frac{6r+1}{\ell} \in \mathbb{Z}_{(\ell)}$, and so, the right hand side reduces to 0. Therefore, we arrive at the desired conclusion, i.e. that
\begin{equation*}
    p_{\alpha}(\ell^2 n + r) \equiv 0 \pmod{\ell^v}.
\end{equation*}

Next, we work out the case of $d= 10$. Similar arguments can be made about $d= 14$ and $26$. Our initial steps are nearly analogous to that of $d= 4$. We once again start by writing $v := \ord_\ell(\aa-10)$ such that $\aa - 10 = \ell^{v} u$ for some $u \in \Z_{(\ell)}$. We also define $r_0$ to be the smallest positive integer such that $12r_0 + 5 \equiv 0 \pmod{\ell}$. We eventually arrive at the analogue of Equation \eqref{eq:3.6}, which is that
\begin{equation} \label{eq:3.7}
     \sum_{n = 0}^\infty p_{\aa}(\ell n+ r_0)q^{12n + \frac{12r_0 +5}{\ell}} \equiv \sum_{n = 0}^\infty a_{10}(\ell n)q^{n}\cdot \poc{q^{12\ell}}^{\ell^{v-2} u}
   \pmod{\ell^{v}}. 
\end{equation}

Recall from Equation \eqref{eq:2.1} that we can write $\eta(12\tau)^{10}$ as linear combinations of two Hecke eigenforms. We have that
\begin{equation*}
    \eta(12\tau)^{10} = \frac{1}{96} ((E_4(12\tau)\eta(12\tau)^2 + 48 \eta(12\tau)^{10}) - (E_4(12\tau)\eta(12\tau)^2 - 48 \eta(12\tau)^{10})).
\end{equation*}

Each of $E_4(12\tau)$, $\frac{\eta(12\tau)^{2}}{q}$, and $\frac{\eta(12\tau)^{10}}{q}$ on the right hand side of Equation \eqref{eq:2.1} are expressions of $q^{4}$. As a result, for $10$-satisfactory primes $\ell$, the $\ell^{\text{th}}$ coefficient in both eigenforms of Equation \eqref{eq:2.1} are $0$. It follows from \Cref{first_lemma} that $a_{10}(\ell k) = 0$ for $k \in \mathbb{Z}_{\ell}$. 

We extract the terms of the form $q^{\ell n + \frac{12r+5}{\ell}}$ from each side of Equation \eqref{eq:3.7}. Once again, because $ord_{\ell}(12r+5) = 1,$ $ord_\ell(\ell n + \frac{12r+5}{\ell})=0$, and so, the right hand side reduces to 0. Thus, we arrive at the desired conclusion that
\begin{equation*}
    p_{\alpha}(\ell^2 n + r) \equiv 0 \pmod{\ell^v}. \qedhere
\end{equation*}
\end{proof}

\begin{proof}[Proof of \Cref{second_theorem}]
The initial steps closely mimic that of the proof of \Cref{first_theorem}. For convenience, we write that $v+1:= \ord_\ell(\aa-2)$ such that $\aa - 2 = \ell^{v+1} u$ for some $u \in \Z_{(\ell)}.$ We relate $p_{\alpha}(n)$ with $\eta(12\tau)^{2}$ through the following steps. Then, we have that
\begin{equation}
    \begin{aligned}
        \sum_{n = 0}^\infty p_{\aa}(n)q^{12n + 1} &= q\poc{q^{12}}^{\alpha} =  q\poc{q^{12}}^{\ell^{v+1} u + 2} \\
    &= q\poc{q^{12}}^{2}\poc{q^{12}}^{\ell^{v+1} u} = \eta(12\tau)^{2}\poc{q^{12}}^{\ell^{v+1} u}.
    \end{aligned}
\end{equation}

Now, applying \Cref{third_lemma} twice, we have that  
\begin{equation}
    \begin{aligned} \label{eq:3.9}
    \sum_{n = 0}^\infty p_{\aa}(n)q^{12n + 1} & = \eta(12\tau)^{2}\poc{q^{12}}^{\ell^{v+1} u} \\
    & \equiv \eta(12\tau)^{2}\poc{q^{12\ell^2}}^{\ell^{v-1} u} \pmod{\ell^v}.
    \end{aligned}
\end{equation}

We rewrite Equation \eqref{eq:3.9} into
\begin{equation}
     \sum_{n = 0}^\infty p_{\aa}(n)q^{12n + 1} \equiv \sum_{n = 0}^\infty a_{2}(n)q^{n}\cdot \poc{q^{12\ell^{2}}}^{\ell^{v-1} u}
   \pmod{\ell^{v}}. 
\end{equation}
Since $\frac{\eta(12\tau)^2}{q}$ is an expression of $q^{12}$, $a_{2}(\ell)= 0$ for $2$-\textit{satisfactory} primes $\ell$. And once again, since  $\eta(12\tau)^2$ is a cuspidal Hecke eigenform, its coefficients are multiplicative among co-prime indices. Therefore, for $k \in \mathbb{Z}_{(\ell)}$, we have that $a_2(\ell k) = 0$. 

Finally, we extract the terms of the form $q^{\ell^2 n + 12r+1}$ from each side. We notice that the right hand side reduces to $0$ as $ord_{\ell}(12r+1)=1$ and arrive at the desired conclusion that
\begin{equation*}
    p_{\alpha}(\ell^2 n + r) \equiv 0 \pmod{\ell^v}. \qedhere
\end{equation*} 
\end{proof}

Before diving into the proof of \Cref{third_theorem}, we prove an auxiliary lemma. 

\begin{lemma}
\label{third_proof_lemma}
    Given a fixed prime $\ell$ and $v \in \mathbb{Z}^{+}$, there exists a $w \in \mathbb{Z}^{+}$ such that $w < \ell^{2v}$ and
    \begin{equation*}
        a_{2}(\ell^w) \equiv 0 \pmod{\ell^v}.
    \end{equation*}
\end{lemma}

\begin{proof} [Proof of \Cref{third_proof_lemma}]
    Because $\frac{\eta(12\tau)^{2}}{q}$ is an expression in terms of $q^{12}$, the statement holds true for $w=1$ when $\ell$ is \textit{$2$-satisfactory}. 
    
    Let $\ell$ be a prime that is not \textit{$2$-satisfactory}. We let $n= \ell^{i}$ for $i \in \mathbb{Z}^{+}$ in \Cref{first_lemma}. Because $\chi(2) = 1$ from \Cref{second_lemma}, it follows that
    \begin{equation} \label{eq:3.11}
        a_{2}(\ell^{i+1}) = a_{2}(\ell^i)a_{2}(\ell) - a_{2}(\ell^{i-1}).
    \end{equation}
    Equation \eqref{eq:3.11} displays a recursion on the sequence of $a_{2}(\ell^i)$ for $i \in \mathbb{Z}^{+} \cup \{0\}$. Notice that the sequence is periodic with respect to modulo $\ell^v$ due to the pigeon hole principle. It follows that the length of the period is at most $\ell^{2v}$, and we let $s \leq \ell^{2v}$ denote the length of the period.
    
    Moreover, it can be observed that the period begins at $a_{2}(1)$. To prove this, assume for the sake of contradiction that the period does not begin at $a_{2}(1)$. We let the first term of the period be $a_2(\ell^{c})$ for some $c > 0$. Then, rearranging Equation \eqref{eq:3.11} and letting $k=c+1$ gives
    \begin{equation}
        \begin{aligned}
            a_{2}(\ell^{c-1}) & \equiv a_{2}(\ell^{c})a_{2}(\ell) - a_{2}(\ell^{c+1}) \\ & \equiv a_{2}(\ell^{c+s})a_{2}(\ell) - a_{2}(\ell^{c+s+1}) \equiv a_{2}(\ell^{c+s-1}) \pmod{\ell^v}. 
        \end{aligned}
    \end{equation}
    This is contradictory to our assumption that $a_2(\ell^c)$ is the first term of the period. Thus, we conclude that the period begins at $a_{2}(1)$.
   
    Now, notice that 
    \begin{equation*}
        a_{2}(\ell^{s-1}) \equiv a_{2}(\ell^{s})a_{2}(\ell)-a_{2}(\ell^{s+1}) \equiv a_{2}(\ell^{0})a_{2}(\ell) - a_{2}(\ell^{1}) \equiv 0 \pmod{\ell^v}.
    \end{equation*}
     As $a_2(1)=1$, setting $w= s-1$ in the statement of the lemma completes the proof. 
\end{proof}

\begin{proof}[Proof of \Cref{third_theorem}] We choose $w$ such that $a_{2}(\ell^w) \equiv 0 \pmod{\ell^v}$, which we know exists by \Cref{third_proof_lemma}. Write $\aa - 2 = \ell^{v+w} u$ for some $u \in \Z_{(\ell)}.$ It follows that 
\vspace{-3mm}
\begin{equation}
    \begin{aligned}
        \sum_{n = 0}^\infty p_{\aa}(n)q^{12n + 1} &= q\poc{q^{12}}^{\alpha} =  q\poc{q^{12}}^{\ell^{v+w} u + 2} \\
    &= q\poc{q^{12}}^{2}\poc{q^{12}}^{\ell^{v+w} u} = \eta(12\tau)^{2}\poc{q^{12}}^{\ell^{v+w} u}. 
    \end{aligned}
\end{equation}

Apply \Cref{third_lemma} $w+1$ times to arrive at    
\begin{equation} \label{eq:3.14}
    \begin{aligned}
        \sum_{n = 0}^\infty p_{\aa}(n)q^{12n + 1} & = \eta(12\tau)^{2}\poc{q^{12}}^{\ell^{w+v} u} \\ & \equiv \eta(12\tau)^{2}\poc{q^{12\ell^{w+1}}}^{\ell^{v-1} u}
     \pmod{\ell^v}. 
    \end{aligned}
\end{equation}

Since $\eta(12\tau)^2$ is a cuspidal Hecke eigenform, we have that $a_2(\ell^w k) \equiv 0 \pmod{\ell^v}$ for all $k \in \mathbb{Z}_{(\ell)}$. 
As $ord_{\ell}(12r+1)=w$, extracting the terms of the form $q^{\ell^{w+1}n + 12r+1}$ from both sides of Equation \eqref{eq:3.14} gives for all $n$ that
\begin{equation*}
    p_{\alpha}(\ell^{w+1}n+r) \equiv 0 \pmod{\ell^v}. \qedhere
\end{equation*}
\end{proof}

\noindent {\bf Acknowledgement.} This project was completed as part of the 2019 Emory REU. The author would like to thank Ken Ono, Larry Rolen, and Ian Wagner for their guidance. The author would also like to thank Erin Bevilacqua, Kapil Chandran, Alice Lin, Eleanor McSpirit, and Junyao Peng for useful discussions. This research was generously supported by the Spirit of Ramanujan Global STEM Talent Search and the Asa Griggs Candler Fund.


\begin{thebibliography}{10}\footnotesize

\bibitem{ahlono} S. Ahlgren and K. Ono, Congruence properties for the partition function, {\it Proc. Natl. Acad. Sci. USA} {\bf 98}, 12882-12884.

\bibitem{atkin} A. O. L. Atkin, Proof of a conjecture of Ramanujan, {\it Glasg. Math. J} {\bf 8}, 14-32.

\bibitem{eigenform} A. Carney, A. Etropolski, and S. Pitman, Powers of the eta-function and hecke operators, {\it Int. J. Number Theory} {\bf 8}, 599-611.

\bibitem{chan} H. H. Chan and L. Wang, Fractional powers of the generating               function for the partition function, {\it Acta Arith.} {\bf 187}, 59-80.

\bibitem{martin} Y. Martin, Multiplicative $\eta$-quotients, {\it Trans. Amer. Math. Soc.} {\bf 348}, 4825-4856.

\bibitem{web} K. Ono, {\it The Web of Modularity: Arithmetic of the Coefficients of Modular Forms and q-series}, American Mathematical Society, Providence, 2004.

\bibitem{lacunary} J. P. Serre, Sur la lacunarit\'e des puissances de $\eta$, {\it Glasg. Math. J} {\bf 27}, 203-221.


\bibitem{watson} G. N. Watson, Ramanujans Vermutung uber Zerfallungsanzahlen, {\it J. Reine Angew. Math} {\bf 179}, 97-128.



\end{thebibliography}
\end{document}